\documentclass[final,leqno]{siamltex}

\usepackage{fullpage}
\usepackage{color}
\usepackage{mathtools}
\usepackage{amssymb}
\usepackage{enumerate}
\usepackage{hyperref}
\usepackage{color, colortbl}
\usepackage{graphicx}
\usepackage{caption}
\usepackage{multirow}
\usepackage{subfigure}
\usepackage{amsmath}

\usepackage{algpseudocode}
\usepackage{algorithm}
\usepackage{footnote}
\usepackage{threeparttable}

\definecolor{Blue}{rgb}{0,0,1}
\definecolor{Red}{rgb}{1,0,0}
\definecolor{Purple}{rgb}{0.5, 0, 0.5}
\definecolor{grey}{gray}{0.6}

\usepackage{soul}

\definecolor{CellHighlighter}{rgb}{1,1,0}
\definecolor{CellBlue}{rgb}{0,1,1}

\newcommand{\hp}{\hat{p}}
\newcommand{\tp}{\tilde{p}}

\def\r{{\vec{r}}}

\newcommand{\hu}{\hat{u}}
\newcommand{\tu}{\tilde{u}}
\newcommand{\bb}{\textbf{b}}

\def\u{{\vec{u}}}
\def\v{{\vec{v}}}

\newcommand{\sA}{\mathcal{A}}
\newcommand{\sR}{\mathcal{R}}

\algnewcommand{\LineComment}[1]{\State \(\triangleright\) #1}

\algnewcommand{\IIf}[1]{\State\algorithmicif\ #1\ \algorithmicthen}
\algnewcommand{\EndIIf}{\unskip\ \algorithmicend\ \algorithmicif}
\algnewcommand{\IFor}[1]{\State\algorithmicfor\ #1\ \algorithmicdo}
\algnewcommand{\EndIFor}{\unskip\ \algorithmicend\ \algorithmicfor}
\newcommand{\RR}[1]{\mathbb{R}^{#1}}
\newcommand{\bz}{\textbf{z}}
\newcommand{\vf}{\textbf{f}}
\def\u{{\vec{u}}}
\def\v{{\vec{v}}}

\def\r{{\vec{r}}}
\newcommand{\spanof}[1]{\text{span}(#1)}

\usepackage{bibentry}
\nobibliography*
\title{Numerical Solution of the Steady-State Navier-Stokes Equations using Empirical Interpolation Methods\thanks{This work was supported by the U.\, S.\, Department of 
Energy Office of Advanced Scientific Computing Research, Applied Mathematics 
program under Award Number DE-SC0009301,
and by the U.\, S.\, National Science Foundation under grant DMS1418754.}}

\author{Howard C. Elman\thanks{Department of Computer Science and Institute for Advanced Computer Studies, University of Maryland, College Park, MD 20742 (elman@cs.umd.edu).} \and Virginia Forstall\thanks{Applied Mathematics and Statistics, and Scientific Computation Program, Department of Mathematics, University of Maryland, College Park, MD 20742 (vhfors@gmail.com).}}

\begin{document}
\maketitle
\begin{abstract}
Reduced-order modeling is an efficient approach for solving parameterized discrete partial differential 
equations when the solution is needed at many parameter values. 
An offline step approximates the solution space and an online step utilizes this approximation, the reduced
basis, to solve a smaller reduced problem at significantly lower cost, producing an accurate estimate of the
solution. 
For nonlinear problems, however, standard methods do not achieve the desired cost savings. 
Empirical interpolation methods represent a modification of this methodology used for cases 
of nonlinear operators or nonaffine parameter dependence. 
These methods identify points in the discretization necessary for representing the nonlinear component of 
the reduced model accurately, and they  
incur online computational costs that are independent of the spatial dimension $N$. 
We will show that empirical interpolation methods can be used to significantly reduce 
the costs of solving parameterized versions of the Navier-Stokes equations, and that iterative 
solution methods can be used in place of direct methods to further reduce the costs of solving
the algebraic systems arising from reduced-order models.
\end{abstract} 
\begin{keywords}reduced basis, empirical interpolation, iterative methods, preconditioning \end{keywords}


\section{Introduction}
Methods of reduced-order modeling are designed to obtain the numerical solution of parameterized partial differential equations (PDEs) efficiently. In settings where solutions of parameterized PDEs are required for many parameters, such as uncertainty quantification, design optimization, and sensitivity analysis, the cost of obtaining high-fidelity solutions at each parameter may be prohibitive. In this scenario, reduced-order models can often be used to keep computation costs low by projecting the model onto a space of smaller dimension with minimal loss of accuracy. 

We begin with a brief statement of the \textit{reduced basis method} for constructing a reduced-order model. Consider an algebraic system of equations $G(u) = 0$ where $u : = u(\xi)$ is an unknown vector of dimension $N$, and $\xi $ is vector of $m$ input parameters. We are interested in the case where this system arises from the discretization of a PDE and $N$ is large, as would be the case for a high-fidelity discretization.  We will refer to this system as the \textit{full model}. We would like to compute solutions for many parameters $\xi$. Reduced basis methods compute a (relatively) small number of full model solutions, $u(\xi_1)$, \ldots $u(\xi_k)$, known as snapshots, 
and then for other parameters, $\xi \ne \xi_j$, construct approximations of $u(\xi)$ in the space spanned by $\{ u(\xi_j) \}_{j=1}^k$. In the offline-online paradigm, the offline step, which may be expensive, computes the snapshots using traditional (PDE) solvers. The offline step builds a basis of the low-dimensional vector space spanned by the snapshots. The online step, which is intended to be inexpensive (because $k$ is small), uses a projected version of the original problem (determined, for example, by a Galerkin projection) in the $k$-dimensional space. The projected problem, known as the \textit{reduced model}, has a solution $\tilde{u}(\xi)$ which is an approximation of the solution $u(\xi)$. 

A straightforward implementation of the reduced-basis method is only possible for linear problems that have affine dependence on the parameters. Such problems have the form $G(u) = 0 $ where
\begin{equation}
G(u) = A(\xi)u - b = \left(\sum_{i=1}^{l} \varphi_i(\xi) A_i \right) u - b
\end{equation}
and $\{A_i\}_{i=1}^l$ are parameter-independent matrices. 
Let $Q$ be a matrix of dimensions $N \times k$ whose columns span the space spanned by the snapshots. For example, $Q$ can be taken to be an orthogonal matrix obtained using the Gram-Schmidt process applied to $[u(\xi_1),..., u(\xi_k)]$; construction of $Q$ is part of the offline step.
With this decomposition, the reduced model obtained from a Galerkin condition is $G^r(\hu) = 0$ where
\begin{equation} 
G^r(\hu) = Q^T A(\xi) Q \hat{u} - Q^T b= \left(\sum_{i=1}^l \varphi_i(\xi) (Q^T A_i Q)\right)\hu - Q^T b \; . 
\label{eq:assembly}
\end{equation}
Computation of the matrices $\{Q^T A_iQ\}$ can be included as part of the offline step. With this preliminary computation, the online step requires only the summation of the terms in equation~(\ref{eq:assembly}), an $O(lk^2)$ operation, and then the solution of the system of order $k$. Clearly, the cost of this online computation is independent of $N$, the dimension of the full model. 

However, when this approach is applied to a nonlinear problem, the reduced model is not independent of the dimension of the full model. 
Consider a problem with a nonlinear component $F(u( \xi))$, so the full model is 
\begin{equation}
G(u(\xi)) = Au(\xi) + F(u(\xi)) - b = 0 \; .
\label{eq:full}
\end{equation}
The reduced model obtained from the Galerkin projection is
\begin{equation}
G^r(\hu(\xi)) =  Q^T A Q \hat{u}(\xi) + Q^T F(Q \hat{u}(\xi)) - Q^T b  = 0 
\label{eq:reduced}
\end{equation}
Although the reduced operator  $Q^T F(Q\hat{u}( \xi))$ is a mapping from $\RR{k} \to \RR{k}$, any nonlinear solution algorithm (e.g.\@ Picard iteration) requires the evaluation of the operator $F(Q\hat{u}(\xi))$ as well as the multiplication by $Q^T$. Both computations have costs that depend on $N$, the dimension of the full model. 


Empirical interpolation methods 
\cite{barrault2004empirical, chaturantabut2010nonlinear, grepl2007efficient,Maday-Nguyen-Patera-Pau}
use interpolation to reduce the cost of the online construction 
for nonlinear operators or nonaffine parameter dependence.  
The premise of these methods is to interpolate the operator 
using a subset of indices from the model. 
The interpolation depends on an empirically derived basis that can also 
be constructed as part of an offline procedure. 
This ensures that $F(Q\hat{u}(\xi))$ is evaluated only at a relatively small 
number of indices. 
These values are used in conjunction with a separate basis constructed to 
approximate the nonlinear operator. 
The efficiency of this approach also depends on the fact that for all 
$i$, $F_i(u(\xi))$ depends on a relatively small, $O(1)$, number of components 
of $u$. 

Computing the solution of the reduced model for a nonlinear operator requires a nonlinear iteration based on a linearization strategy, which requires the solution of a reduced linear system at each step. Thus, each iteration has two primary costs, the computation of the Jacobian corresponding to $Q^T F(u(\xi))$, and the solution of the linear system at each step of the nonlinear iteration. 
Empirical interpolation addresses the first cost, by using an approximation of 
$Q^T F (u(\xi))$. 
To address the second cost, one option is to use direct methods to solve the reduced linear systems. 
In \cite{elman2015preconditioning}, however, we have seen that iterative methods are effective for 
solving reduced models of linear operators of a certain size. 
In this paper, we extend this approach, using preconditioners that are precomputed in the offline stage, to nonlinear problems solved using empirical interpolation.  
We explore this approach using a Picard iteration for the linearization strategy. 

We will demonstrate the efficiency of combining empirical interpolation with an iterative 
linear solver by computing solutions of the steady-state incompressible Navier-Stokes equations with 
random viscosity coefficient:
 \begin{equation}
 \label{eq:ns}
  \renewcommand{\baselinestretch}{1}
\begin{array}{rclcc}
- \nabla \cdot \nu(\cdot,\xi)\nabla \u(\cdot,\xi) + \u(\cdot,\xi) \cdot \nabla \u(\cdot,\xi) + \nabla p(\cdot,\xi) &=& f  & \text{ in } & D \times \Gamma \\  \nonumber
\nabla \cdot \u(\cdot,\xi) &=& 0 & \text{ in } & D \times \Gamma \\ 
\u(\cdot,\xi) & =& b & \text{ on }& \partial  D \times \Gamma \; ,  \nonumber
\end{array}
\end{equation}
where $\u(\cdot, \xi)$ is the flow velocity, $p(\cdot,\xi)$ is the scalar pressure, and $b$ determines the Dirichlet boundary conditions, and the viscosity coefficient satisfies $\nu(\cdot,\xi) > 0$. Models of this type have been used to model the viscosity in multiphase flows \cite{kim2007phase, olshanskii2006analysis, tan2008immersed}. The boundary data $b$ and forcing function $f$ could also be parameter-dependent, although we will not consider such examples here. 
 
An outline of this paper is as follows. 
In Section \ref{sec:deim}, we outline the details of the empirical interpolation strategy 
that we use, the so-called discrete empirical interpolation method (DEIM) 
\cite{chaturantabut2010nonlinear}.
(See \cite[Ch.\ 10]{Quarteroni-etal-RB} for discussion and comparison of different 
variants of this idea.)
In Section \ref{sec:ns} we introduce the steady-state Navier-Stokes equations with an 
uncertain viscosity coefficient and describe the full, reduced, and DEIM models for 
this problem. 
We present numerical results in Section \ref{sec:ex} including a comparison of snapshot 
selection methods for DEIM and a discussion of accuracy of the DEIM. 
In addition, we discuss a generalization of this approach known as a gappy-POD method 
\cite{carlberg2013gnat, everson1995karhunen}. 
Finally, in Section \ref{sec:iterative_deim}, we discuss the use of iterative methods 
for solving the reduced 
linear systems that arise from DEIM. 
This includes a presentation of new preconditioning techniques for use in this setting 
and discussion of their effectiveness.

\section{The discrete empirical interpolation method}
\label{sec:deim}
The discrete empirical interpolation method utilizes an approximation $\bar{F}(u)$ of a nonlinear function $F(u)$ \cite{chaturantabut2010nonlinear}. The keys to efficiency in this algorithm are that
\begin{enumerate}
\item only a small number of indices of $F$ are used in each component 
\item each component of the nonlinear function depends only on a few indices of the input variable. 
\end{enumerate} 
The key to the accuracy of this method is to select the indices of the discrete PDE that are most important to produce an accurate representation of the nonlinear component of the solution projected on the reduced space. The efficiency requirement is clearly satisfied when the nonlinear function is a PDE discretized using the finite element method \cite{antil2013application}.

Given a PDE that depends on a set of parameters $\xi = [\xi_1, ..., \xi_m]^T$ with solution $u(\cdot,\xi)$, the \textit{full model} is the discretized equation such that $G(u(\cdot,\xi)) = 0$.  The offline step for the traditional approach to the reduced basis method takes full solutions at several parameters and constructs a matrix $Q$ of rank $k$ such that $\text{range}(Q) = \text{span}\{u(\cdot,\xi^{(1)}), ..., u(\cdot,\xi^{(k)})\}$, where the solutions $u(\cdot,\xi^{(1)}), ..., u(\cdot,\xi^{(k)})$ are known as \textit{snapshots}. 
The online step approximates the true solution $u$ with an approximation $\tu \approx Q \hu$. Using the Galerkin projection, the reduced model, $G^r$, is
$$G^r(\hu) = Q^T G(Q \hat{u}) \; . $$ and 
the Jacobian of $G^r(\hu)$ is
$$J_{G^r}(\hu) \ = \  \frac{\partial G^{r}(\hu)}{\partial \hu} \ = \  
Q^T \frac{\partial}{\partial \hu} G(\tu(\hu))\ = \  
Q^T \frac{\partial G }{\partial \tu }\frac{\partial \tu }{\partial \hu} \ = \ 
Q^T J_G(\tu) Q =  Q^T J_G(Q\hu) Q \; . $$ 

As observed above, when the operator is linear and affinely dependent on the parameters, the online costs, (of forming and solving the reduced system (\ref{eq:assembly})) are independent of $N$. This is not true for nonlinear or nonaffine operators.
Consider Newton's method for the reduced model with a nonlinear operator: 
\begin{equation}
\hu_{n+1} = \hu_{n} - J_{G^r}(\hu_n)^{-1} G^r(\hu_n) 
\label{eq:iteration}
\end{equation}
Each iteration in equation (\ref{eq:iteration}) requires the construction of $J_{G^r}(\hu_n)= Q^T J_G(Q\hu_n)Q$. The construction of the matrix $J_G(Q\hu_n)$ as well as the multiplications by $Q$ and $Q^T$ have costs that  depend on $N$.

In the DEIM approach, the snapshots $u(\xi^{(1)}), ..., u(\xi^{(k)})$ are obtained from the full model 
and the reduced basis is constructed to span these snapshots. 
In addition, DEIM requires a separate basis to represent the nonlinear component of 
the solution. This basis is constructed using a matrix of snapshots of function values 
$S = [F(u(\xi^{(1)})), F(u(\xi^{(2)})), ..., F(u(\xi^{(s)}))]$.
Then, using methods similar to finding the reduced basis, a basis is chosen to approximately span the space spanned by these snapshots. One approach for doing this is to use a proper orthogonal decomposition (POD) of the snapshot matrix $S$
\begin{equation}
 S = \bar{V} \Sigma W^T 
 \label{eq:svd_intro} 
\end{equation}
where the singular values in $\Sigma$ are sorted in order of decreasing magnitude and 
$\bar{V}$ and $W$ are orthogonal. 
DEIM will use columns of $\bar{V}$ to approximate $F(\tilde{u}(\xi))$. 
It may happen that $n_{deim}<k$ columns of $\bar{V}$ are used here, giving 
a submatrix $V$ of $\bar{V}$.
We will discuss the method used to select $n_{deim}$ in Section \ref{sec:ex}.

Given the nonlinear basis from $V$, the DEIM selects indices of $F$ so the interpolated nonlinear component on the range of $V$ in some sense represents a good approximation to the complete set of values of $F(u(\xi))$. In particular, the approximation of the nonlinear operator is
$$\bar{F} (u(\xi))  = V(P^T V)^{-1}P^T F(u(\xi))$$
where $P^T$ extracts entries of $F(u)$ corresponding to the interpolation points from the spatial 
grid.\footnote{An implementation does not literally construct the matrix $P$; instead, 
an index list is used to extract the required entries of $V$.} 
This approximation satisfies $P^T \bar{F} = P^T F$. To construct $P$, a greedy procedure is used to minimize the error compared with the full representation of $F(u)$ \cite[Algorithm 1]{chaturantabut2010nonlinear}. 
For each column of $V$, $v_i$, the DEIM algorithm selects the row index for which the difference between the column $v_i$ and the approximation of $v_i$ obtained using the DEIM model with nonlinear basis and 
the first $i-1$ columns of V is maximal, that is,
the index of the maximal entry of $r = v_i - \hat{V}(P^T \hat{V})^{-1} P^T v_i$, where $\hat{V}$ denotes the first $i-1$ columns of $V$. We present this in Algorithm \ref{alg:deim}. 

\begin{algorithm}
\renewcommand{\baselinestretch}{1}
\caption{DEIM  \cite{chaturantabut2010nonlinear}}
Input: $V = [v_1, ..., v_{n_{deim}}]$, an $N \times n_{deim}$ matrix with columns made up of the left singular vectors from the POD of the nonlinear snapshot matrix $S$. \\
Output: $P$, extracts the indices used for the interpolation.
\begin{algorithmic}[1]
\State $\rho = \text{argmax}(|v_1|)$, the index of the maximal entry of $|v_1|$ \\
 $\widehat{V} =  [v_1]$, $P = [e_{\rho}]$ 
\For {$i =2:n_{deim}$} \\
\hspace{\algorithmicindent} Solve $(P^T \widehat{V})c = P^T v_i$ for $c$  \\
\hspace{\algorithmicindent} $r = v_i - \widehat{V}c$			\\
 \hspace{\algorithmicindent} $\rho = \text{argmax}(|r|)$		\\
 \hspace{\algorithmicindent} $\widehat{V} = [\widehat{V}, v_i]$, $P = [P,e_{\rho}] $			
 \EndFor 
\end{algorithmic}
\label{alg:deim}
\end{algorithm}

Incorporating this approximation into the reduced model, equation~(\ref{eq:reduced}), yields
\begin{equation}
\bar{F}^r(\hu) =  Q^T \bar{F}(\tu)  = Q^T V ( P^T V) ^{-1}  P^T F(Q\hu) \; . 
\label{reduced_deim}
\end{equation}
The construction of nonlinear basis matrix $V$ and the interpolation points are part of the offline computation. Since $L^T = Q^T V(P^T V)^{-1}$ is parameter independent, it too can be computed offline. Therefore, the online computations required are to compute $P^T J_F(u)$ and assemble $L^T (P^T J_F(u))Q$. For $P^TJ_F(u)$, we need only to compute the components of $J_F(u)$ that are nonzero at the interpolation points. This is where the assumption that each component of $F(u)$ (and thus $J_F(u)$) depends on only a few entries of $u$ is utilized. With a finite element discretization, a component $F_i(u)$ depends on the components $u_j$ for which the intersection of the support of the basis functions have measure that is nonzero. See \cite{antil2013application} for additional discussion of this point. The elements that must be tracked in the DEIM computations are referred to as the sample mesh. When the sample mesh is small, the computational cost of assembling $L^T (P^T J_F(u)) Q$ scales not with $N$ but with the number of interpolation points. Therefore, DEIM will decrease the online cost associated with assembling the nonlinear component of the solution. 

For the Navier-Stokes equations, the nonlinear component is a function of the velocity. We will discretize the velocity space using biquadratic ($Q_2$) elements. In this case, an entry in $F_i(u)$ depends on at most nine entries of $u$. Thus this nonlinearity is amenable to using DEIM. An existing finite element routine can be used for the assembly of the required entries
of the Jacobian using the sample mesh, a subset of the original mesh, as the input. 

The accuracy of this approximation is determined primarily by the quality of the nonlinear 
basis $V$. 
This can be seen by considering the error bound  
\begin{equation}
 || F - \bar{F} ||_2 \le || (P^T V )^{-1} ||_2 ||(I - VV^T)F||_2
 \label{eq:error_bound}
 \end{equation}
which is derived and discussed in more detail in \cite[Section 3.2]{chaturantabut2010nonlinear}. There it is shown that the greedy selection of indices in Algorithm \ref{alg:deim} limits the growth of $||(P^T V)^{-1}||_2$ as the dimension of $V$ grows. The second term $||(I - VV^T)F||_2$ is the quantity that is determined by the quality of $V$. Note that if $V$ is taken from the truncated POD of $S$, the matrix of nonlinear snapshots, then $||(I-VV^T)S||_F^2$ is minimized \cite{antil2013application} where $||\cdot||_F$ is the Frobenius norm ($||X||_F^2  = \sum_i \sum_j |x_{ij}|^2$). So the accuracy of the DEIM approximation depends on two factors. The first  is the number, $n_{deim}$, of singular vectors kept in the POD. The truncated matrix $V \Sigma_{deim}W_{deim}^T$ is the optimal rank-$n_{deim}$ approximation of $S$, but a higher rank approximation will improve accuracy of the DEIM model. In fact the error $||(I - VV^T)F||_2$ approaches 0 in the limit as $n_{deim}$ approaches $N$. The second factor is the quality of the nonlinear snapshots in $S$. The nonlinear component should be sampled well enough to capture the variations of the nonlinear component throughout the solution space. 
A comparison of methods for selecting the snapshot set is given in Section \ref{sec:offline}. 

Given the full model defined in equation (\ref{eq:full}), let $J_G(u) = A + J_F(u)$  
denote the Jacobian matrix. The Jacobian of the reduced model equation (\ref{eq:reduced}) is then
$$J_{G^r}(\hu) = Q^T J_G(Q\hu) Q = Q^T A Q + Q^T J_F(Q\hu) Q, $$
and the Jacobian $\bar{F}^r(\hu)$ of (\ref{reduced_deim}) is
$$J_{\bar{F}^r}(\hu) \ = \  Q^T J_{\bar{F}}(Q\hu) Q  \  = \ Q^T V (P^T V)^{-1} P^T J_F(u) Q \; . $$


\section{Steady-state Navier-Stokes equations}
\label{sec:ns}
A discrete formulation of the steady-state Navier-Stokes equations~(\ref{eq:ns}) is to find $\u_h \in X_E^h$ and $p_h \in M^h$ such that
\begin{equation}
 \renewcommand{\baselinestretch}{1}
\begin{array}{rclc}
(\nu(\cdot,\xi) \nabla \u_h, \nabla \v_h) + (\u_h \cdot \nabla \u_h,\v_h) - (p_h,\nabla \cdot \v_h) &=& (f,\v_h) &\forall \v_h \in X_0^h \\ \nonumber
(\nabla \cdot \u_h,q_h) &=& 0 &\forall q_h \in M^h
\end{array}
\end{equation}
where $X_E^h$ and $M^h$ are finite-dimensional subspaces of the Sobolev spaces $H_0^1=\{\vec{v}\in H^1|\vec{v} = 0 \text{ on } \partial D\}$ and $L_2(D)$; see \cite{elman2014finite, Girault-Raviart} for details. We will use div-stable $Q_2$-$P_{-1}$ finite element (biquadratic velocities, piecewise constant discontinuous pressure). Let $\{\phi_1, ..., \phi_{n_u}\}$ represent a basis of $Q_2$ and $\{\psi_1, ..., \psi_{n_p}\}$ represent a basis of $P_{-1}$. 

We define the following vectors and matrices, where $u$
and $p$ here represent the vectors of coefficients that determine $\u_h$ and $p_h$,
respectively:
\begin{equation*}
\bz = \begin{bmatrix} u \\ p \end{bmatrix},
\end{equation*}
\vspace{-.12in}
\begin{equation*}
[A(\xi)]_{ij} = \int \nu(\xi) \nabla  \phi_i : \nabla \phi_j \,,   \; \; \; \; \;  
[B]_{ij} = - \int \psi_i ( \nabla \cdot \phi_j )\,, \; \; \; \; \; 
[N(u)]_{ij} = \int (\u_h \cdot \nabla \phi_j) \cdot \phi_i\,,
\end{equation*}
\vspace{-.12in}
\begin{equation*}
[\vf]_i = (f, \phi_i), \; \; \; \; \;  
[g(\u)]_i =  - (\nabla \cdot \u_h, \psi_i) ,  \; \; \; \; \; 
\bb(\xi) = \begin{bmatrix} \vf - A(\xi)u_{bc} \\ g(u_{bc}) \end{bmatrix}\,,
\end{equation*}
where $u_{bc}$ is the vector of coefficients of the discrete velocity field $\u_{bc}$ that 
interpolates the Dirichlet boundary data $b(\cdot,\xi)$ and is zero everywhere on the interior of the mesh. 
We denote the velocity solution on the interior of the mesh, $\u_{in}$, so that $\u = \u_{bc} + \u_{in}$ and 
$\u_{in}$ satisfies homogenous Dirichlet boundary conditions. 
The reduced basis is constructed using snapshots of $\u_{in}$ so the approximation of the velocity 
solution generated by the reduced model is of the form $\tu = u_{bc} + Q_u \hat{u}$ where the columns 
of  $Q_u$ correspond to a basis spanning the space of velocity snapshots with homogeneous Dirichlet
boundary conditions. 

\subsection{Full model}

With this notation, the full discrete model for the Navier-Stokes problem 
with parameter $\xi$ is to find $\bz(\xi)$ such 
that $G(\bz(\xi)) = 0$ where 
\begin{equation}
\renewcommand{\baselinestretch}{1}
G(\bz(\xi)) = \begin{bmatrix}
A(\xi) & B^T \\
B & 0 \\
\end{bmatrix}
\begin{bmatrix}
 u \\  p
\end{bmatrix} 
+ 
\begin{bmatrix}
 N(u) & 0 \\
0 & 0 \\
\end{bmatrix}
\begin{bmatrix}
 u \\  p
\end{bmatrix}
-
\begin{bmatrix}
\vf \\ 0
\end{bmatrix}
\; . 
\end{equation}
We utilize a Picard iteration to solve the full model, monitoring the norm of the nonlinear residual $G(\bz_n(\xi))$ for convergence. The nonlinear Picard iteration to solve this model is described in Algorithm \ref{alg:full}. 
\begin{algorithm}
\caption{Picard iteration for solving the discrete steady-state Navier-Stokes equations}
\begin{algorithmic}[1]
\State The nonlinear iteration is initialized with the solution to a Stokes problem  
\begin{equation}\begin{bmatrix} A(\xi) & B^T \\
B & 0 \\
\end{bmatrix}    \begin{bmatrix} \u_{in,0} \\ p_0 \end{bmatrix}= \bb(\xi) \; .         
\label{eq:full_stokes}
\end{equation}
\item Incorporate the boundary conditions
$$\u_{0} = \u_{bc} + \u_{in,0} \; . $$
\item Solve 
\begin{equation}
\left(  \begin{bmatrix} A(\xi) & B^T \\
B & 0 \\
\end{bmatrix} + \begin{bmatrix} N(\u_n) & 0 \\
0 & 0 \\
\end{bmatrix}  \right) \begin{bmatrix} \delta \u \\ \delta p \end{bmatrix}= -G(\bz_n) \; . 
\label{eq:full_picard}
\end{equation}

\item Update the solutions
\begin{align*}
\u_{n+1} &= \u_{n} + \delta \u \\ \nonumber
p_{n+1} &= p_{n} + \delta p  \; . 
\end{align*}

\item Exit if $$||G(\bz_{n+1})||_2  < \delta \left | \left| \bb(\xi)  \right | \right |_2 \; , $$
otherwise return to step 2. 
\end{algorithmic}
\label{alg:full}
\end{algorithm}

\begin{algorithm}
\renewcommand{\baselinestretch}{1}
\caption{Picard iteration for solving the reduced steady-state Navier-Stokes equations}
\begin{algorithmic}[1]
\State Initialize the Picard iteration by solving the reduced Stokes problem
$$ Q^T \begin{bmatrix} A(\xi) & B^T \\
B & 0 \\
\end{bmatrix} Q   \begin{bmatrix} \hu_{0} \\ \hp_0 \end{bmatrix}= Q^T \bb(\xi) \; .  $$
\item Solve the reduced problem for the Picard iteration
$$
\left(Q^T  \begin{bmatrix} A(\xi) & B^T \\B & 0 \\ \end{bmatrix} Q 
+ Q^T  \begin{bmatrix} N(\tu_n) & 0 \\ 0 & 0 \\ \end{bmatrix}Q  \right) 
\begin{bmatrix} \delta \hu \\ \delta \hp \end{bmatrix}
= -G^r(\tilde{\bz}_n) \; . 
$$
Note that the when the dependence on the parameters is affine, the first term in the left hand 
side can be computed primarily offline as in equation (\ref{eq:assembly}). 
\item Update the reduced solutions
\begin{align*}
\hu_{n+1} &= \hu_{n} + \delta \hu\\
\hp_{n+1} &= \hp_{n} + \delta \hp \; . \\
\end{align*}
\item Update the approximation to the full solution
\begin{align*}
\tu_{n+1} &= \u_{bc} + Q_u \hu_{n+1} \\ 
\tp_{n+1} &= Q_p \hp_{n+1}  \; . \\
\end{align*}
\item Compute $N(\tu_{n+1})$. 
\item Compute $G^r(\tilde{\bz}_{n+1})$. 
\item Exit if
$$||G^r(\tilde{\bz}_{n+1})||_2 < \delta  \left| \left| Q^T \bb(\xi) \right| \right|_2 \; , $$
otherwise return to step 2.
\end{algorithmic}
\label{alg:reduced}
\end{algorithm}


\subsection{Reduced model} 
Next, we present a reduced model that does not make use of the DEIM strategy. This is not meant to be practical strategy, but it is presented as a comparison to illustrate how the additional approximation used for DEIM affects both the accuracy of solutions obtained using DEIM and the speed with which they are obtained. Offline we compute a reduced basis 
\begin{equation*}
\renewcommand{\baselinestretch}{1}
Q = \begin{bmatrix} Q_u & 0 \\ 0 & Q_p  \end{bmatrix} \; ,
\end{equation*}
where $Q_u$ represents the reduced basis of the velocity space and $Q_p$ the reduced basis for the 
pressure space. 
We defer the details of this offline construction to Section \ref{sec:offline}.  
For a given $Q$, the Galerkin reduced model is
$$
G^r(\bz) \ = \ Q^T G(\bz)  
\ = \ 
\left( Q^T  \begin{bmatrix}
A(\xi)  & B^T \\
B & 0 \\
\end{bmatrix}
Q \right )
\begin{bmatrix}
 \hu \\  \hp
\end{bmatrix} 
+ 
\left (Q^T 
\begin{bmatrix}
 N(\tu) & 0 \\
0 & 0 \\
\end{bmatrix}
Q \right )
\begin{bmatrix}
\hu \\  \hp
\end{bmatrix}
-
Q^T 
\begin{bmatrix}
\vf \\ 0
\end{bmatrix}
\; . 
$$
Using the nonlinear Picard iteration, the reduced model is described in Algorithm \ref{alg:reduced}. 

After the convergence of the Picard iteration determined by $G^r(\tilde{\bz}_{n+1})$, we compute the ``full'' residual $G(\tilde{\bz}_n)$. Note that this residual is computed only once; because the cost of computing it is $O(N)$, it is not monitored during the course of the iteration. The full residual indicates how well the reduced model approximates the full solution, so it is used to measure the quality of the reduced model via the error indicator
\begin{equation}
\eta_{\xi} = ||G(\tilde{\bz}_n(\xi))||_2/\left|\left| \bb(\xi) \right|\right|_2 \; . 
\label{eq:residual}
\end{equation}

\subsection{DEIM model}
The DEIM model has the structure of the reduced model but with the nonlinear component $F$ replaced by the approximation $\bar{F}$. First in the offline step, we compute $V$, $P$, and $L^T =  Q_u^T V ( P^T V) ^{-1}$. The DEIM model is
\begin{equation}
\renewcommand{\baselinestretch}{1}
G^{deim}(\bz) = 
\left( Q^T \begin{bmatrix}
A(\xi) & B^T \\
B & 0 \\
\end{bmatrix}
Q \right )
\begin{bmatrix}
 \hu \\  \hp
\end{bmatrix} 
+ 
\begin{bmatrix}
L^TP^T N(\tu)Q_u & 0 \\
0 & 0 \\
\end{bmatrix}
\begin{bmatrix}
 \hu \\  \hp
\end{bmatrix}
- Q^T 
\begin{bmatrix}
\vf \\ 0
\end{bmatrix}
\; . 
\end{equation}
The computations are shown in Algorithm \ref{alg:deim_model}. 
Recall from the earlier discussion of DEIM that $P^T N(u)$ is not computed by forming 
the matrix $N(u)$. 
Instead, only the components of $N(u)$ corresponding to the indices required by 
$P^T$ are constructed, using the elements of the discretization mesh that contribute to 
those indices.
Note that error indicator $\eta_{\xi}$ in equation (\ref{eq:residual}) depends on $G(\tilde{\bz})$. 
This quantity contains $N(\tu_n)$ and not $P^T N(\tu_n)$, so that 
computing it requires assembly of $N(\tu_n)$ on the entire mesh.
As for the reduced model, to avoid this expense, this computation is performed only after 
convergence of the nonlinear iteration. 

\begin{algorithm}
\renewcommand{\baselinestretch}{1}
\caption{DEIM model for the steady-state Navier-Stokes equations}
\begin{algorithmic}[1]
\State Initialize the Picard iteration by solving the reduced Stokes problem
\begin{equation}
Q^T \begin{bmatrix} A(\xi) & B^T \\
B & 0 \\
\end{bmatrix} Q   \begin{bmatrix} \hu_{0} \\ \hp_0 \end{bmatrix}= Q^T \bb(\xi)  \; . 
\label{eq:deim_stokes}
\end{equation}
\item Solve the reduced problem for the Picard iteration
\begin{equation}
\left(Q^T  \begin{bmatrix} A(\xi) & B^T \\
B & 0 \\
\end{bmatrix} Q + \begin{bmatrix}  L^T (P^T N(\tu_n))Q_u & 0 \\
0 & 0 \\
\end{bmatrix}  \right) \begin{bmatrix} \delta \hu \\ \delta \hp \end{bmatrix}= -G^{deim}(\tilde{\bz}_n) \; . 
\label{eq:deim_picard}
\end{equation}
Note that the term on the left can be computed cheaply as in equation~(\ref{eq:assembly}) so we only need to update the upper left corner of the matrix as the Picard iteration proceeds. \label{line:linear_solve}
\item Update the reduced solutions
\begin{align*}
\hu_{n+1} &= \hu_{n} + \delta \hu\\
\hp_{n+1} &= \hp_{n} + \delta \hp \; . \\
\end{align*}
\item Compute  $P^TN(\tu_{n+1})= P^TN(Q\hu_{n+1})$ at required indices.  
\item Compute $G^{deim}(\tilde{\bz}_{n+1})$. 
\item Exit if
$$||G^{deim}(\tilde{\bz}_{n+1})||_2 < \delta  \left| \left|Q^T \bb(\xi) \right| \right|_2 \; , $$
otherwise, return to step 2.
\end{algorithmic}
\label{alg:deim_model}
\end{algorithm}


\subsection{Inf-sup condition}
\label{sec:infsup}
We turn now to the construction of the reduced basis 
\begin{equation*}
\renewcommand{\baselinestretch}{1}
Q = \begin{bmatrix} Q_u & 0 \\ 0 & Q_p  \end{bmatrix} \; .
\end{equation*}
Given $k_s$ snapshots of the full model, a natural choice is to have the following spaces generated by these snapshots, 
\begin{equation}
\renewcommand{\baselinestretch}{1}
\label{eq:space}
\begin{array}{rcl}
\text{span}(Q_u) &=& \text{span} \{ \u_{in}(\xi^{(1)}), ...,\u_{in}(\xi^{(k_s)}) \} \\ 
\text{span}(Q_p) &=& \text{span} \{ p(\xi^{(1)}), ...,p(\xi^{(k_s)}) \}  \; . \\ 
\end{array}
\end{equation}
However, this choice of basis does not satisfy an inf-sup condition
\begin{equation*}
\gamma_R := \min_{0\ne q_R \in \text{span}(Q_p)} \max_{0\ne\v_R \in \text{span}(Q_u)} 
\frac{(q_R,\nabla \cdot \v_R)}{|\v_R|_1 ||q_R||_0} \ge \gamma^* > 0
\end{equation*}
where $\gamma^*$ is independent of $Q_u$ and $Q_p$. 
To address this issue, we follow the enrichment procedure of 
\cite{quarteroni2007numerical,rozza2007on}. 
For $i = 1, ..., k_s$, let $\r_h(\cdot,\xi^{(i)})$ be the solution to the Poisson problem 
 \begin{equation*}
 (\nabla\r_h(\cdot,\xi^{(i)}),\nabla\v_h) = (p_h(\cdot,\xi^{(i)}),\nabla\cdot\v_h)  \; \; \forall \v_h \in X_0^h \; , 
 \end{equation*}
 and let $Q_u$ of equation (\ref{eq:space}) be augmented by the corresponding discrete solutions 
 $\{\r(\xi^{(i)})\}$, giving the enriched space
  $$\text{span}(Q_u) = \text{span} \{ \u_{in}(\xi^{(1)}), ...,\u_{in}(\xi^{(k_s)}),\r(\xi^{(1)}), ..., \r(\xi^{(k_s)} ) \} \; . $$
These enriching functions satisfy 
\begin{equation*}
\r_h(\cdot,\xi^{(i)}) = \arg \sup_{\v_h \in X_0^h} \frac{(p_h(\cdot,\xi^{(i)}), 
\nabla \cdot \v_h)}{|\v_h|_1}  \; , 
\end{equation*}
and thus $\gamma_R$ defined for the enriched velocity space, $\spanof{Q_u}$, together with $\spanof{Q_p}$, satisfies the inf-sup condition
\begin{equation*}
\gamma_R \ge \gamma_h := \min_{0\ne q_h \in M_h} \max_{0\ne\v_h \in X_0^h} \frac{(q_h,\nabla \cdot \v_h)}{|\v_h|_1 ||q_h||_0} \; . 
\end{equation*}

\section{Experiments}
\label{sec:ex}

We consider the steady-state Navier-Stokes equations (\ref{eq:ns}) for driven cavity flow posed on a square domain $D = (-1,1) \times (-1,1)$. The lid, the top boundary ($y = 1$), has velocity profile
$$u_x  = 1 - x^4, \; \;  u_y = 0 \; , $$
and no-slip conditions $\u = (0,0)^T$ hold on other boundaries. The source term is $f \equiv 0$.  The discretization is done on a 
div-stable $n \times n$ $Q_2$-$P_{-1}$ 
(biquadratic velocities, discontinuous piecewise constant pressures)
element grid, giving a discrete velocity space of order $(n+1)^2$ and pressure space of order $3(n/2)^2$. 

To define the uncertain viscosity, we divide the domain $D$ into $m = n_d \times n_d$ subdomains as seen in Figure~\ref{fig:domain}, and the viscosity is taken to be constant and random on each subdomain, $\nu(\xi) = \xi_i$. The random parameter vector, $\xi = [\xi_1, ..., \xi_m]^T \in \Gamma$, is comprised of uniform random variables such that $\xi_i \in \Gamma_i = [0.01,1]$ for each $i$. Therefore, the local subdomain-dependent Reynolds number, $\sR = 2/\nu$, will vary between 2 and 200 for this problem. Constant Reynolds numbers in this range give rise to stable steady solutions \cite{elman2014finite}.

\setlength{\unitlength}{1cm}
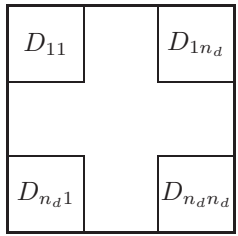
\begin{figure}[ht]
\begin{center}
\begin{picture}(3,3)
\put(0,2){\framebox(1,1)[c]{$D_{11}$}}
\linethickness{.25pt}
\put(2,2){\framebox(1,1)[c]{$D_{1n_d}$}}
\linethickness{.5pt}
\put(2,0){\framebox(1,1)[c]{$D_{n_d n_d}$}}
\put(0,0){\framebox(1,1)[c]{$D_{n_d1}$}}
\put(0,0){\framebox(3,3)[s]{}}
\end{picture}
\end{center}
\caption{Flow domain with piecewise random coefficients for viscosity.}
\label{fig:domain}
\end{figure}

The implementation uses IFISS \cite{ifiss} to generate the finite element matrices for the full model. The matrices are then imported into Python and the full, reduced, and DEIM models are constructed and solved using a Python implementation run on an Intel 2.7 GHz i5 processor and 8 GB of RAM. The full model is solved with the method described in equation~(\ref{eq:sing_fix}) using sparse direct methods implemented in  the UMFPACK suite \cite{davis2004algorithm} for the system solves in equations~(\ref{eq:full_stokes}) and~(\ref{eq:full_picard}). For this benchmark problem (with enclosed flow), these linear systems are singular \cite{bochev2005finite}. This issue is addressed by augmenting the matrix, for example that of (\ref{eq:full_stokes}), as
\begin{equation}
\begin{bmatrix} A(\xi) & B^T & 0 \\
B & 0 & 0 \\
0 & \bar{p} & 0 \\
\end{bmatrix}  \begin{bmatrix} \u_{in} \\ p \\ z \end{bmatrix}= \begin{bmatrix} \multirow{2}{*}{\bb($\xi$)} \\
  \\0  \end{bmatrix}  \; , 
\label{eq:sing_fix}
\end{equation}
where $\bar{p}$ is a vector corresponding to the element areas of the pressure elements. This removes the singularity by adding  a constraint via a Lagrange multiplier so that the average pressure of the solution is zero \cite{silvester2013ifiss}. 
The same constraint is added to the systems in equation~(\ref{eq:full_picard}). 

\subsection{Construction of $Q$ and $V$}
\label{sec:offline}
\begin{algorithm}
\renewcommand{\baselinestretch}{1}
\caption{Construction of reduced basis $Q$ via random sampling, construction of nonlinear basis $V$}
Cost: $n_{trial}$ reduced problems and $k$ full problems. 
\begin{algorithmic}[1] 
\LineComment Compute the reduced basis and the nonlinear snapshots
\State Solve the full problem $G(\bz(\xi^{(0)}))=0$ to tolerance $\delta$ for $\bz_n(\xi^{(0)})$.\\
Compute the enriched velocity, $\vec{r}(\xi^{(0)})$.\\
Initialize $Q_u = [\u_{in,n}(\xi^{(0)}), \vec{r}(\xi^{(0)})]$ and $Q_p = [p_n(\xi^{(0)})]$. \\
Save the nonlinear component of the solution $S = [N(\u_n)\u_n]$. \\
$k_s=0$
\For {$i =1:n_{trial}$} \\
\hspace{\algorithmicindent}Randomly select $\xi^{(i)}$.  \\
\hspace{\algorithmicindent}Solve reduced model $G^r(\tilde{\bz}(\xi^{(i)})) 
= 0$ to tolerance $\delta$ and compute the residual indicator 
$\eta_{\xi^{(i)}}$.
\hspace{\algorithmicindent} \If  {$\eta_{\xi^{(i)}} > \tau$}\\
\hspace{\algorithmicindent}\hspace{\algorithmicindent} $k_s = k_s + 1$\\
\hspace{\algorithmicindent}\hspace{\algorithmicindent} 
Solve full model $G(\bz(\xi^{(i)})) = 0$. \\
\hspace{\algorithmicindent}\hspace{\algorithmicindent} 
Compute the enriched velocity, $\vec{r}(\xi^{(i)})$. \\
\hspace{\algorithmicindent}\hspace{\algorithmicindent} 
Add $\u_{in,n}(\xi^{(i)})$ and  $\vec{r}(\xi^{(i)})$ to $Q_u$ and 
$p_n(\xi^{(i)})$ to $Q_p$ using modified Gram-Schmidt.\\
\hspace{\algorithmicindent}\hspace{\algorithmicindent} 
Add the nonlinear component to the matrix of nonlinear snapshots, 
$S = [S, N(\u_n)\u_n]$. 
\EndIf
\EndFor
\LineComment The nonlinear snapshot matrix is $S = [F(\bz(\xi^{(1)})), ..., F(\bz(\xi^{(k_s)}))]$.
\label{line:nonlinear_basis}
\State Compute the POD of the nonlinear snapshot matrix: $$S = \bar{V}\begin{bmatrix}\sigma_{1} && \\ &\ddots&\\ &&\sigma_{k_s}\end{bmatrix}W^T . $$ \\
Choose $n_{deim}=rank(S)$.\\
Define $V = \bar{V}[:,1:n_{deim}]$. 
\State Compute $P$ using Algorithm \ref{alg:deim} with input $V$.  \label{line:call_deim}
\end{algorithmic}
\label{alg:adaptive_offline}
\end{algorithm}

We now describe the methodology we used to compute the reduced bases, $Q_u$ and $Q_p$, and the nonlinear basis $V$. The description of the construction of $Q_u$ and $Q_p$ is presented in Algorithm~\ref{alg:adaptive_offline}. The reduced bases $Q_u$ and $Q_p$ are constructed using random sampling of $n_{trial}$ samples of $\Gamma$, denoted $\Gamma_{trial}$. The bases are constructed so that all samples $\xi \in \Gamma_{trial}$ have a residual indicator, $\eta_{\xi}$, less than a tolerance, $\tau$. The procedure begins with single snapshot $\bz(\xi^{(0)})$ where $\xi^{(0)} = E(\xi)$. The bases are initialized using this snapshot, such that $Q_u= [\u_{in,n}(\xi^{(0)}), \r(\xi^{(0)})]$  and $Q_p = [p_n(\xi^{(0)})]$. Then for each sample of $\Gamma_{trial}$, the reduced-order model is solved with the current bases $Q_u$ and $Q_p$. The quality of the reduced solution produced by this reduced-order model can be evaluated using the error indicator $\eta_{\xi}$ of (\ref{eq:residual}). If $\eta_{\xi}$ is smaller than the tolerance $\tau$, the computation proceeds to the next sample. If the error indicator exceeds the tolerance, the full model is solved, and then the new snapshots, $u_{in,n}(\xi)$ and $p_n(\xi)$, and the enriched velocity $\r(\xi)$, are used to augment $Q_u$ and $Q_p$. The experiments use $\tau = 10^{-4}$ and $n_{trial} = 2000$ parameters to produce bases $Q_u$ and $Q_p$. 

An alternative to this strategy of random sampling is greedy sampling, which 
produces a basis of quasi-optimal dimension \cite{binev2011convergence, buffa2012priori}, in the 
sense that the the maximum error differs from the Kolmogorov n-width by a constant factor.
Our experience \cite{elman2015preconditioning} is that the performance of the random sampling strategy 
used here is comparable to that of a greedy strategy.  
In particular, for several linear benchmark problems, to achieve comparable 
accuracy, we found that the size of the reduced basis generated by sampling 
was never more than 10\% larger than that produced by a greedy algorithm and 
in many cases the basis sizes were identical. 
The computational cost (in CPU time) of the random sampling strategy is 
significantly lower. 
(See a discussion of this point in Section \ref{sec:online_component}.)
Since our concern in this study is online strategies for reducing the cost 
of the reduced model, we use the random sampling strategy for the offline 
computation and remark that the online solution strategies considered here
can be used for a reduced basis obtained using any method.

We turn now to the methodology for determining the nonlinear basis $V$, which is the truncated form of 
$\bar{V}$ defined in equation~(\ref{eq:svd_intro}). 
It was shown in \cite{chaturantabut2010nonlinear} that the choice of the nonlinear basis, $V$
in equation (\ref{eq:error_bound}), is important for the accuracy of the DEIM model.
DEIM uses a POD approach for constructing the nonlinear basis. 
This POD has a two inputs, $S$, the matrix of nonlinear snapshots, and $n_{deim}$, the number of vectors retained after truncation. 
We will compare three strategies for sampling $S$ specified in 
line \ref{line:nonlinear_basis} of Algorithm \ref{alg:adaptive_offline}.
\begin{enumerate}
\item \textbf{Full}($n_{trial}$). This method is most similar to the method used to generate the nonlinear basis in \cite{chaturantabut2010nonlinear}. The matrix of nonlinear snapshots, $S$, is computed from the full solution at every random sample (i.e. $\{N(\u(\xi^{(i)})) \u(\xi^{(i)})\}_{i=1}^{n_{trial}}$). 
This computation is part of the offline step but its cost, for solving  solving $n_{trial}$ full 
problems, can be quite high.
\item \textbf{Full}($k_s$). This is the sampling strategy included in Algorithm~\ref{alg:adaptive_offline}. It saves the nonlinear component only when the full model is solved for augmenting the reduced basis, $Q$. Therefore the snapshot set $S$ contains $k_s$ snapshots. 
\item \textbf{Mixed}($n_{trial})$. The final approach aims to mimic the Full($n_{trial}$) method with less offline work. This method generates a nonlinear snapshot for each of the $n_{trial}$ random samples using full solutions when they are available (from full solution used for augmenting the reduced basis) and reduced solutions when they are not. As the reduced basis is constructed, when the solution to the full problem is not needed for the reduced basis (i.e. when $\eta_{\xi^{(i)}} < \tau$), use the reduced solution $\tu(\xi^{(i)})$ to generate the nonlinear snapshot $N(\tu(\xi^{(i)})) \tu(\xi^{(i)})$, where $\tu(\xi^{(i)}) = \u_{bc} + Q_u \hat{u}(\xi^{(i)})$ and $Q_u$ is the basis at this value of $i$. 
Thus $S$ contains $n_{trial}$ snapshots, but it is constructed using only $k_s$ full model solves.
\end{enumerate}

\begin{figure}[hbt]
	\centering
	\subfigure{\includegraphics[width=0.7\textwidth]{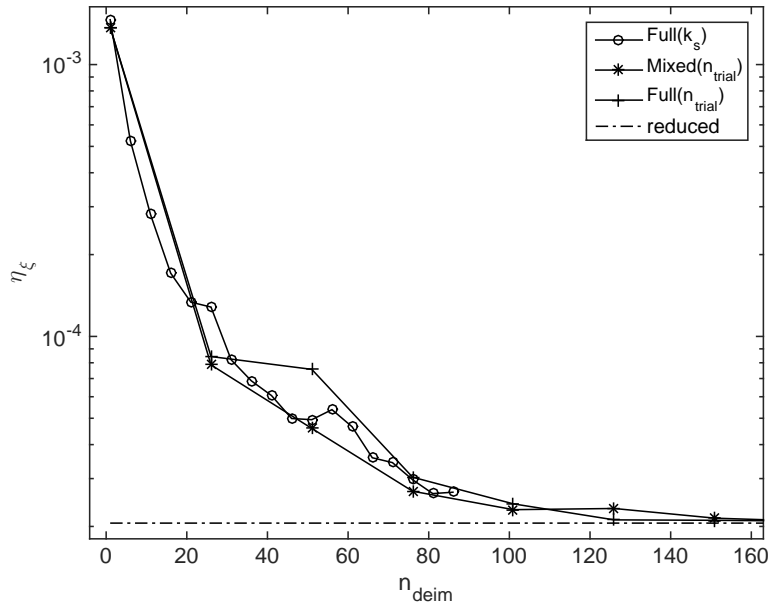}}
       \caption{A comparison of methods to generate nonlinear snapshots for the DEIM method. DEIM residual
       vs.~$n_{deim}$ averaged for $n_s = 100$ samples, $n = 32$, $m = 4$, $\tau = 10^{-4}$, $k =306 $, 
       $n_{deim}$ varies. $Q_u$, $Q_p$ are generated using the Algorithm~\ref{alg:adaptive_offline}.  }
\label{fig:offline_comp}
\end{figure}

Figure \ref{fig:offline_comp} compares the performance of the three methods for generating $S$ when Algorithm \ref{alg:adaptive_offline} is used to generate $Q_u$ and $Q_p$. For each $S$, we take the SVD and truncate with a varying number of vectors, $n_{deim}$, and plot the average of the residuals of the DEIM solution for 100 randomly generated samples of $\xi$. The average residual for the reduced model without DEIM is also shown. 
It can be seen that as $n_{deim}$ increases, the residuals of the DEIM models approach the residual 
obtained without using DEIM.
It is also evident that the three methods perform similarly. 
Thus, the Mixed($n_{trial}$) approach provides accurate nonlinear snapshots with fewer full solutions than
the Full($n_{trial}$) method.
The Full($k_s$) method is essentially as effective as the others and requires both fewer
full-system solves than Full($n_{trial}$) and the SVD of a smaller matrix than Mixed($n_{trial}$);
we used Full($k_s$) for the remainder of this study.
For this method, $n_{deim}\le k_s$ (in fact, in all cases $n_{deim}=k_s$), which is why its 
results do not fully extend across the horizontal axis in Figure \ref{fig:offline_comp};
if the (slightly higher) accuracy exhibited by the other methods is needed, it could be obtained 
using the Mixed($n_{trial}$)  method at relatively little extra cost.

\subsection{Online component - DEIM model versus reduced model}
\label{sec:online_component}
In Section~\ref{sec:deim}, we presented analytic bounds for how accurately the DEIM approximates the nonlinear component of the model. To examine how the approximation affects the accuracy of the reduced model, we will compare the error indicators of DEIM with those obtained using the reduced model without DEIM. 
We perform the following offline and online computations:
\begin{enumerate}
\item 
Offline: Use Algorithm \ref{alg:adaptive_offline} with input $\tau = 10^{-4}$ and $n_{trial} = 2000$
to generate the reduced bases $Q_u$, $Q_p$, $V$ and the indices in $P$. 
\item Online: Solve the problem using the full model, reduced model without DEIM, and the reduced model 
with DEIM for $n_s=10$ sample parameter sets. 
\end{enumerate}
Table \ref{tab:redvdeim} presents the results of these tests.
For each benchmark problem, there are four entries.
Three of the entries, the first, third and fourth, show the time required to 
meet the stopping criterion in each
of Algorithms 2 (full system), 3 (reduced system) and 4 (reduced system
with DEIM) with tolerance $\delta =10^{-8}$, 
together with the final relative residual of the full nonlinear system.
The reduced models use a reduced residual for the stopping test; as we
have observed, this is done for DEIM to make the cost of the iteration
depend on the size of the reduced model rather than $N$, the size of the 
full model.  
(Although the reduced model without DEIM has costs that depend on $N$, we 
used the same stopping criterion in order to assess the difference between
the two reduced models.)
The other entries of Table \ref{tab:redvdeim},  the second for each example,
are the time and residual data for the full system solve with a milder tolerance,
$\delta = 10^{-4}$; with this choice, the relative residual is comparable in size to that obtained
for the reduced models.

For each test problem, the smallest computational time is in boldface.
The times presented for the reduced and DEIM models are the CPU time spent 
in the online computation for the nonlinear iteration only and do not include 
assembly time or the time to compute the full nonlinear residual after the 
iteration for the reduced model has converged. 

The results demonstrate the tradeoff between accuracy and time for the 
three models. 
For example, in the case of $m = 16$ and $n = 65$, the spatial dimension and 
parameter dimension are large enough that the DEIM model is fastest, and 
this is generally true for higher-resolution models.
If smaller residuals are needed, the accuracy of the DEIM model can be improved
by either increasing $n_{deim}$ or improving the accuracy of the reduced model. 
Increasing $n_{deim}$ has negligible effect on online costs (since $L^T$ is 
computed offline),  but the benefits of doing this are limited.
As can be seen from Figure 4.2, the accuracy of the DEIM solution is limited by the accuracy 
of the reduced solution;
once that accuracy is obtained, increasing $n_{deim}$ will produce no additional improvement.
The accuracy of the reduced model can be improved using a stricter 
tolerance $\tau$ during the offline computation. 
Since choosing a stricter tolerance for the reduced model causes the size 
of the reduced basis k to increase, the cost of the reduced and DEIM models 
will also increase.
This means that the benefits of more highly accurate DEIM computations will be obtained only for 
higher-resolution models.

In Table \ref{tab:redvdeim}, the results for $n = 32$ and $m \ge 25$ are not shown. For these problems, the number of snapshots required to construct the reduced model, $k_s$ exceeds the size of the pressure space $n_p = 3(n/2)^2 = 768$. This means that the number of snapshots required for the accuracy of the velocity is higher than the number of degrees of freedom in the full discretized pressure space. Therefore, the spatial discretization is not fine enough for reduced-order modeling to be necessary.  

{\bf Remark.}  Although offline computations are typically viewed as being inconsequential, we 
also note that for some of the larger examples tested, where $n = 128$, the costs of the 
offline construction are significant.
First, the full solutions require over 2 minutes of CPU time. For $m = 36$, 1013 full 
solutions and 2000 reduced solutions were required. 
The cost of each full solve is 132 seconds and the costs of the reduced solves are 
as high as 98.1 seconds. 
The computation time for the assembly of the reduced models are not presented in 
Table \ref{tab:redvdeim} and are also high. 
Since $Q$ is changing during the offline stage, the assembly process cannot be made 
independent of $N$.  
For these experiments, the offline computation for $m=49$ took approximately 
five days, in contrast to 25 seconds for the DEIM solution on the finest grid.


\begin{table}[htp]
\renewcommand{\baselinestretch}{1}
\resizebox{\textwidth}{!}{
\begin{tabular}{|c|l|cc|cc|cc|cc|cc|cc|}
\hline
$n$& \multicolumn{1} {c|}{$m$}   & \multicolumn{2}{c|}{4} & \multicolumn{2}{c|}{9}& \multicolumn{2}{c|}{16} & \multicolumn{2}{c|}{25}&\multicolumn{2}{c|}{36}&\multicolumn{2}{c|}{49} \\
\hline
\multirow{6}{*}{$32$}   & \multicolumn{1} {c|}{$k$}             &  \multicolumn{2} {c|}{306}    &\multicolumn{2} {c|}{942}& \multicolumn{2} {c|}{1485} &  \multicolumn{2} {c|}{} &  \multicolumn{2} {c|}{}
& \multicolumn{2} {c|}{}   \\
& \multicolumn{1} {c|}{$n_{deim}$}              &       \multicolumn{2} {c|}{102}&\multicolumn{2} {c|}{314}     &\multicolumn{2} {c|}{495}      &\multicolumn{2} {c|}{} &\multicolumn{2} {c|}{}         &\multicolumn{2} {c|}{}   \\
\cline{2-14}
&& time & res  & time & res & time & res  & time & res  & time & res  & time & res \\
\cline{2-14}
&\multirow{2}{*}{Full}   &1.11   & $<$1.E-08      & 1.16 &        $<$1.E-08        & 1.05 &       $<$1.E-08 & & & && &      \\
&          &0.44   & 3.78E-05      & 0.52 &        2.51E-05        & \textbf{0.66} &       8.89E-06 & & & && &      \\
&Reduced&0.15   & 2.22E-05      & 1.10 &        3.18E-05        & 3.40 &        4.61E-05 & & & & & &    \\
&DEIM   &\textbf{0.06}  & 2.46E-05      &\textbf{ 0.51} &       3.41E-05        & 1.86 &        4.75E-05 &&&&&& \\
\hline

\multirow{6}{*}{$64$}   & \multicolumn{1} {c|}{$k$}             &  \multicolumn{2} {c|}{273}    &\multicolumn{2} {c|}{825}& \multicolumn{2} {c|}{1503   }
&\multicolumn{2} {c|}{2394} &  \multicolumn{2} {c|}{3339 } & \multicolumn{2} {c|}{4455}  \\
& \multicolumn{1} {c|}{$n_{deim}$}              &\multicolumn{2} {c|}{91}       &\multicolumn{2} {c|}{275}      &\multicolumn{2} {c|}{501}  &\multicolumn{2} {c|}{798} &\multicolumn{2} {c|}{1113} &\multicolumn{2} {c|}{1485}  \\
\cline{2-14}
&& time & res & time & res  & time & res &   time & res  & time & res &time & res \\
\cline{2-14}
&\multirow{2}{*}{Full} & 11.0 & $<$1.E-08 &  10.8 &    $<$1.E-08 &  10.2     & $<$1.E-08    &  11.3  & $<$1.E-08 &  10.1 & $<$1.E-08 & 10.3   & $<$1.E-08 \\
& & 4.58 & 4.42E-05 &  5.07 &       2.99E-05  &  5.61       & 1.07E-05       &  \textbf{5.75} & 1.06E-05 &  \textbf{5.70} & 1.44E-05 & \textbf{5.42}         & 3.97E-05 \\
&Reduced & 0.36& 4.94E-05&2.38& 1.67E-05        & 7.13 &4.49E-05& 19.4  &4.90E-05& 37.2 &5.89E-05 &     70.8&7.42E-05    \\
&DEIM   &\textbf{0.07} & 4.55E-05       & \textbf{0.39} &1.69E-05& \textbf{1.57}&       4.57E-05& 6.23& 5.00E-05& 13.4& 5.98E-05& 29.6   &      7.51E-05 \\
\hline
\multirow{6}{*}{$128$}          & \multicolumn{1} {c|}{$k$}             &  \multicolumn{2} {c|}{237}    &\multicolumn{2} {c|}{732}& \multicolumn{2} {c|}{1383} &  \multicolumn{2} {c|}{2109} &  \multicolumn{2} {c|}{3039}&\multicolumn{2} {c|}{4083}   \\
& \multicolumn{1} {c|}{$n_{deim}$}              &\multicolumn{2} {c|}{79}       &\multicolumn{2} {c|}{244}      &\multicolumn{2} {c|}{461}      &\multicolumn{2} {c|}{703}&\multicolumn{2} {c|}{1013}&\multicolumn{2} {c|}{1361} \\
\cline{2-14}
&& time & res & time & res  & time & res &   time & res  & time & res&time&res \\
\cline{2-14}
&\multirow{2}{*}{Full} & 135 & $<$1.E-08 & 141 & $<$1.E-08     & 147 & $<$1.E-08 & 155 &   $<$1.E-08      & 132 & $<$1.E-08 & 148 & $<$1.E-08 \\
& & 55.1 & 3.24E-05 & 54.4 & 4.33E-05       & 53.6 & 4.50E-05& 63.5 &       3.57E-05         & 66.8 & 3.35E-05 & 57.1 &             5.74E-05 \\
&Reduced & 1.14 & 1.38E-05&7.83 &       2.28E-05& 22.2  &2.76E-05& 51.0 &5.17E-05 & 101 &5.33E-05       & 196 &7.05E-05 \\
&DEIM&\textbf{0.11} & 1.41E-05&\textbf{ 0.51}&  2.47E-05        & \textbf{1.76}& 2.80E-05       & \textbf{4.84}& 5.27E-05& \textbf{11.0}&       5.49E-05&       \textbf{24.8} & 7.16E-05 \\
\hline
\end{tabular}}
\caption{Solution time and accuracy for Full, Reduced and DEIM models.}
\label{tab:redvdeim}
\end{table}

\begin{figure}[htb]
	\centering
	\subfigure{\includegraphics[width=0.6\textwidth]{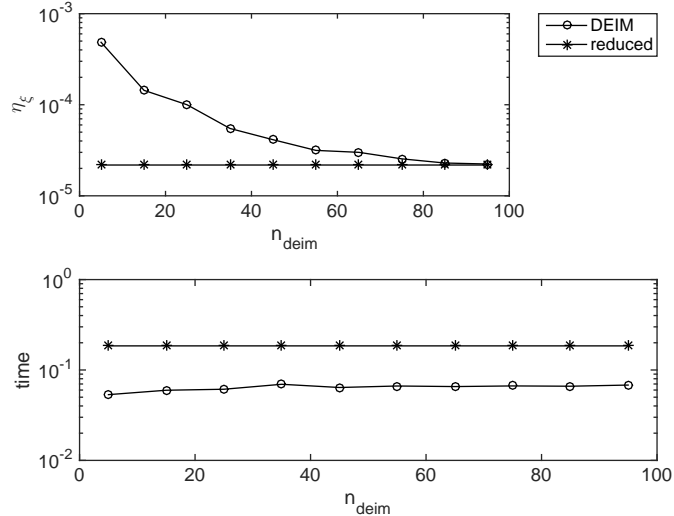}}
       \caption{Top: Error indicator for DEIM model versus $n_{deim}$. Bottom: CPU time to solve using DEIM direct versus $n_{deim}$. For $n = 32$, $m = 4$, $\tau = 10^{-4}$, $k = 306$, averaged over $n_s = 10$ samples.}
       \label{fig:deim_comp}
\end{figure}

Figure~\ref{fig:deim_comp} illustrates the tradeoff between accuracy and time for the DEIM. The top plot compares the error indicators for an average $n_s = 10$ parameters and the bottom plot shows the CPU time for the two methods. While the cost of the DEIM does increase with the number of vectors $n_{deim}$, we reach similar accuracy as for the reduced model at a much lower cost. It can also be seen that the cost of increasing $n_{deim}$ is small since the maximum considered here ($n_{deim}=96$) is significantly smaller than $N = 1089$. 

\subsubsection{Gappy POD}
Another way to increase the accuracy of the reduced model is to increase the number 
of interpolation points in the approximation, while keeping the number of basis vectors 
fixed. 
This alternative to DEIM for selecting the indices is the so-called gappy POD method 
\cite{everson1995karhunen}. 
This method allows the number or rows selected by $P^T$ to exceed the number of columns of $V$. 

The approximation of the function using gappy POD looks similar to DEIM, but the inverse of $P^T V$ 
is replaced with the Moore-Penrose pseudoinverse $(P^T V)^{\dagger}$ \cite{carlberg2013gnat}
\begin{equation} 
\hat{F}(u) = V (P^T V)^{\dagger} P^T F(u) \; .
\end{equation}
To use this operator, we compute $P^T F(u)$ and solve the least squares problem
$$\alpha = \arg\min_{\hat{\alpha}} ||P^T V\hat{\alpha} - P^T F(u) ||_2 \; , $$
which leads to the approximation $\hat{F}(u) = V \alpha$. Like $(P^T V)^{-1}$, the 
pseudoinverse can be precomputed, in this case using the SVD $P^TV = U\Sigma W^T$, giving
$$(P^T V)^{\dagger}=W\Sigma^{\dagger}U^T $$
where $\Sigma^{\dagger}$ is \cite{660book} 
$$\begin{bmatrix}
1/\sigma_{1} && \\ &\ddots && 0\\ &&1/\sigma_{n_{deim}} & \end{bmatrix} \; . 
$$

With this approximation, the index selection method is described in 
Algorithm \ref{alg:gappy}. 
Given $V$ computed as in Algorithm \ref{alg:adaptive_offline},
we require a method to determine the selection of the row indices that will lead to an accurate representation of the nonlinear component. The approach in \cite{carlberg2013gnat} is an extension of the greedy algorithm used for DEIM (Algorithm~\ref{alg:deim}). It takes as an input the number of grid points and the basis vectors. It simply chooses additional indices per basis vector where the indices correspond to the maximum index of the difference of the basis vector and its projection via the gappy POD model. 
Recall that in DEIM, the index associated with vector $v_i$ that maximizes
$|v_i - V(P^TV)^{-1}P^T v_i|$ (called $\rho$ in Algorithm 1) is found, and $P$ is augmented by $e_{\rho}$. 
In this extension of the algorithm, at each step $i$, this type of construction is done $n_g/n_v$ times 
for $v_i$, where an index $\arg\max(|v_i - V(P^TV)^{\dagger}P^T v_i|)$ is found and then P and the 
projection of $v_i$ are updated each time a new such index is chosen. 

\begin{algorithm}
\caption{Index selection using gappy POD \cite{carlberg2013gnat}}
Input: $n_{g}$ number of indices to choose, $V = [v_1, ..., v_{n_{v}}]$, an $N \times n_{v}$ matrix with columns made up of the left singular vectors from the POD of the nonlinear snapshot matrix $S$. \\
Output: $P$, extracts the indices used for the interpolation.
\begin{algorithmic}[1]
\State $n_b = 1$, $n_{it} = \min(n_v,n_{g})$
	$n_{c,\min} = \left \lfloor{\frac{n_v}{n_{it}}}\right \rfloor$
	$n_{a,\min} = \left \lfloor{\frac{n_{g}}{n_v}}\right \rfloor$
	\For {$i =1, ..., n_{it}$} \\
\hspace{\algorithmicindent}$n_c = n_{c,\min}$, $n_a = n_{a,\min}$ 
\hspace{\algorithmicindent} \IIf{$i <=(n_v\mod n_{it}) $} {$n_c = n_c + 1$} \EndIIf
\hspace{\algorithmicindent} \IIf {$i <=(n_{g}\mod n_v)$} $n_a = n_a + 1$\EndIIf
		\If {$i == 1$} \\
\hspace{\algorithmicindent}\hspace{\algorithmicindent}$r = \sum_{q = 1}^{n_c} v_{q}^2$
			\For {$j = 1, ..., n_a$}
				$\rho_j =  \text{argmax}(r)$,
				$r[\rho_j] = 0$ \EndFor \\ 
\hspace{\algorithmicindent}\hspace{\algorithmicindent}$P = [e_{\rho_1}, ..., e_{\rho_{n_a}}]$,
			$\widehat{V} = [v_1, ..., v_{n_c}]$
		\Else
			\For{$q = 1, ..., n_c$}\\
\hspace{\algorithmicindent}\hspace{\algorithmicindent}\hspace{\algorithmicindent}$\alpha = \min_{\hat{\alpha}} ||P^T \widehat{V}\hat{\alpha} - P^T v_{n_b+q}||_2$\\
\hspace{\algorithmicindent}\hspace{\algorithmicindent}\hspace{\algorithmicindent}$R_q  = v_{n_b+q} - \widehat{V}\alpha$ 
			\EndFor \\
\hspace{\algorithmicindent}\hspace{\algorithmicindent}$r = \sum_{q = 1}^{n_c} R_{q}^2$
			\For{$j=1, .., n_a$}\\
\hspace{\algorithmicindent}\hspace{\algorithmicindent}\hspace{\algorithmicindent}$\rho_j =  \text{argmax}(r)$,
		$P = [P, e_{\rho_j}]$
				\For {$q =1, ..., n_c$}\\
\hspace{\algorithmicindent}\hspace{\algorithmicindent}\hspace{\algorithmicindent} \hspace{\algorithmicindent}$\alpha = \min_{\hat{\alpha}}  ||P^T \widehat{V}\hat{\alpha} - P^T v_{n_b+q}||
					_2$\\
\hspace{\algorithmicindent}\hspace{\algorithmicindent}\hspace{\algorithmicindent} \hspace{\algorithmicindent}$R_q = v_{n_b+q} - \widehat{V}\alpha$		
					\EndFor\\
\hspace{\algorithmicindent}\hspace{\algorithmicindent}\hspace{\algorithmicindent}$r = \sum_{q = 1}^{n_c} R_{q}^2$
			\EndFor\\
\hspace{\algorithmicindent}\hspace{\algorithmicindent}$\widehat{V} = [\widehat{V}, v_{n_b+1}, ..., v_{n_b+n_c}]$, $n_b = n_b + n_c$
		\EndIf
		
 \EndFor 
\end{algorithmic}
\label{alg:gappy}
\end{algorithm}

\begin{figure}[hbt]
	\centering
	\subfigure{\includegraphics[width=0.6\textwidth]{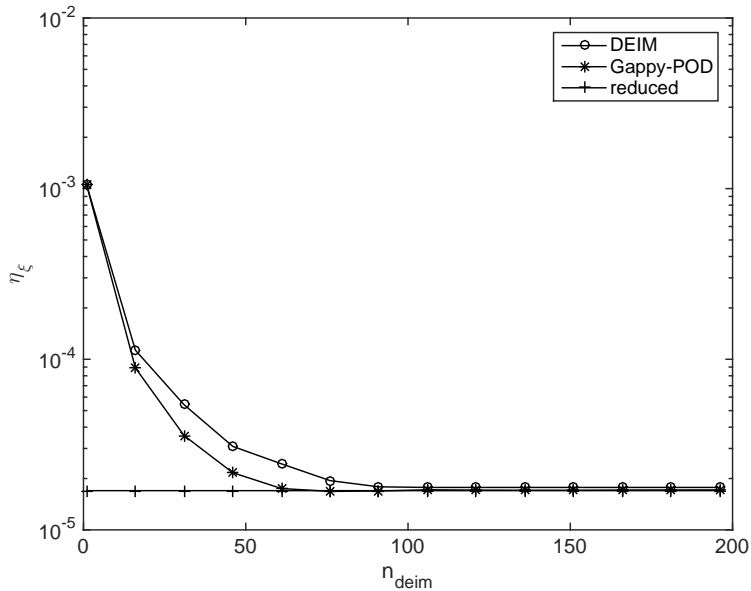}}
       	\caption{Average error indicator as a function of basis vectors for reduced, DEIM, and gappy POD methods. For $n = 32$, $m = 4$, $\tau = 10^{-4}$, $k = 306$, $n_{deim}$ varies, and $n_{g} = 2n_{deim}$. Averaged over $n_s = 100$ samples.}
       \label{fig:gappy_comp4}
\end{figure}

To compare the accuracy of this method with DEIM, we use Algorithm \ref{alg:adaptive_offline} to 
compute DEIM and modify line~\ref{line:call_deim} to use Algorithm \ref{alg:gappy} with 
$n_g= 2n_{deim}$ for a range of values of $n_{deim}$, 
so that there are two indices selected for each $v_i$ in the gappy algorithm. 
We use both methods to approximate the nonlinear component and solve the resulting models. 
Figure \ref{fig:gappy_comp4} shows the error indicators as functions of $n_{deim}$ for 
both methods.
It is evident that for smaller number of basis vectors the gappy POD provides additional accuracy. 
However, for larger number of basis vectors, the additional accuracy provided by gappy POD is small. 
Thus, the gappy POD method can be used to improve the accuracy when the number of basis vectors is limited. 
Since $S$ is generated using the Full($k_s$) described in Section~\ref{sec:offline}, no additional 
accuracy is gained for the DEIM method when $n_{deim} > 102$. 

\section{Iterative methods} 
\label{sec:iterative_deim}
We have seen that the DEIM generates reduced-order models that produce solutions 
essentially as accurate as the reduced solution. 
In addition, Table \ref{tab:redvdeim} and Figure \ref{fig:deim_comp} illustrate that as expected, 
the DEIM model significantly decreases the online time spent constructing the nonlinear component 
of the reduced model. 
Since $Q_u^T N(\u) Q_u$ has been replaced by a cheap approximation, $L^T P^T N(\u)Q_u$, 
the remaining cost of the nonlinear iteration in the DEIM is that of the linear system solve 
in line \ref{line:linear_solve} of Algorithm \ref{alg:deim_model}. 

The cost of this computation depends on the rank of the reduced basis $k$, 
which is $3k_s$ in our setting, where $k_s$ is the number of snapshots used 
to construct the reduced basis. 
These quantities depend on properties of the problem such as the number of 
parameters in the model or the desired level of accuracy in the reduced model.
In contrast, the cost of solving the full model is independent of the number of parameters, and  
it could be as small as $O(N)$ if multigrid methods can be utilized.
Thus, it may happen that $k$  is much less than $N$, but the $O(k^3)$ cost of solving the
DEIM model  is larger than the cost of solving the full model. 
This is the case, for example, for  $n=64$, $m = 49$ in Table \ref{tab:redvdeim}, where 
the CPU time to solve the full model is half that of solving the DEIM model using 
direct methods. 



An alternative is to use iterative methods to solve the reduced model.
Their cost is $O(k^2p)$ where $p$ is the number of iterations required 
for convergence, so that there are values of $k$ where, if $p$ is small enough,
 iterative methods will be preferable to direct methods. 
 We have seen examples of this for linear problems in \cite{elman2015preconditioning}.
 In this section, we discuss the use of iterative methods based on preconditioned Krylov 
 subspace methods to improve the efficiency of the DEIM model. 

For iterative methods to be efficient, effective preconditioners are needed. 
In the offline-online paradigm, it is also desirable to make the construction of the
preconditioner independent of parameters, so that this construction can be part of the offline
step.
Thus, we will develop preconditioners that depend only on the mean value $\xi^{(0)}$ from the 
parameter set, and refer to such techniques as ``offline'' preconditioners.
To see the impact of this choice, 
we will also compare the performance of these approaches with ``online'' versions of them, where the 
preconditioning operator for a model with parameter $\xi$ is built using that parameter;
this approach is not meant to be used in practice since its online cost will depend on $N$, but it provides
insight concerning a lower bound on the iteration count that can be achieved using an offline preconditioner.  
Versions of the offline approach have been used with stochastic Galerkin methods in \cite{powell2012preconditioning}.

We consider two preconditioners of the DEIM model: 
\begin{enumerate}
\item \textit{The Stokes preconditioner} is the matrix used for the reduced Stokes solve in equation (\ref{eq:deim_stokes}), 
\begin{equation}
\renewcommand{\baselinestretch}{1}
M_r(\xi) = Q^T \begin{bmatrix} 
A(\xi) & B^T \\ 
B  & 0 \\
\end{bmatrix}
Q \; . 
\end{equation}
\item \textit{The Navier-Stokes preconditioner} uses the converged solution of the full model $\u_n$ as 
the input for $P^TN(\u_n)$ and uses the operator from the DEIM model from 
equation (\ref{eq:deim_picard}),
\begin{equation}
\renewcommand{\baselinestretch}{1}
M_r(\xi) = Q^T \begin{bmatrix} 
A(\xi) & B^T \\ 
B  & 0 \\
\end{bmatrix}
Q  +  \begin{bmatrix} 
L^T P^T N(\u_n) Q_u & 0 \\ 
0  & 0 \\
\end{bmatrix} \; . 
\end{equation}
\end{enumerate}
In both cases, the preconditioned coefficient matrix is $\mathcal{A}(\xi) M_r^{-1}$ where 
$\mathcal{A}(\xi)$ is the coefficient matrix in (\ref{eq:deim_picard}).
According to the comments above, we consider two versions of each of these, the offline version, $M_r(\xi_0)$, and the online version, $M_r(\xi)$. Clearly, for the first step of the reduced nonlinear iteration, a Stokes solve, only one linear iteration will be required using the online Stokes preconditioner.


For these experiments, we solve the steady-state Navier-Stokes equations for the driven cavity 
flow problem using the full model, reduced model, and the DEIM model. 
For the DEIM model, the linear systems are solved using both direct and iterative methods. 
The iterative methods use the preconditioned \texttt{bicgstab} method \cite{Vorst-book}. 

The offline construction is described in Algorithm~\ref{alg:adaptive_offline}; we use $\tau  = 10^{-4}$ and $n_{trial} = 2000$. The algorithm chooses $k_s$ snapshots and produces $Q_u$ of rank $2k_s$ and $Q_p$ of rank $k_s$, yielding reduced models of rank $k = 3k_s$. 
The online experiments are run for $n_s = 10$ random parameters. The average number of iterations required for the convergence of the linear systems is presented in Table~\ref{tab:it_countstab} and the average time for the entire nonlinear solve of each model is presented in Table~\ref{tab:timingstab}. The nonlinear solve time includes the time to compute $N$ or $P^T N$, but not the time for assembly of the linear component of the model nor the computation time for $\eta_{\xi}$, the norm of the full residual (\ref{eq:residual}) for the approximate solution found by the reduced and DEIM models. The iterative methods presented in this table use offline preconditioners.  The lowest online CPU time is in bold. The nonlinear iterations are run to tolerance $\delta = 10^{-8}$ and the \texttt{bicgstab} method for a reduced system with coefficient matrix, $Q^T \sA Q$, stops when the solution $x^{(i)}$ satisfies
$$ \frac{||r - Q^T \sA Qx^{(i)}||}{||r||} < 10^{-9} \;  $$
where $r$ is $Q^T \begin{bmatrix}\mathbf{f}\\ 0 \end{bmatrix}$. 

\begin{table}
\renewcommand{\baselinestretch}{1}
\begin{center}
\begin{tabular}{|c|l|cccccc|}
\hline
$n$& \multicolumn{1} {c}{$m$}	 & 4 & 9 & 16 & 25 & 36 & 49 \\ 
\hline
\multirow{6}{*}{$32$}  & \multicolumn{1} {c|}{$k$} &306 & 942 & 1485 &&& \\
		& \multicolumn{1} {c|}{$n_{deim}$} &102 & 314 &  495 &&& \\
\cline{2-8}
&Offline Stokes 	&11.3&16.2&19.8&&&\\
&Online Stokes	&2.0&2.3&2.3&&&\\
&Offline Navier-Stokes			&11.5&16.1&20.5&&&	\\
&Online Navier-Stokes			&1.8&1.9&1.9&&&		\\	
\hline	
\multirow{6}{*}{$64$} &\multicolumn{1}{c|}{$k$} &273&825& 1503&2394&3339 &4455\\
  & \multicolumn{1} {c|}{$n_{deim}$}	 &91	&275& 501&798 &  1113 &1485 \\
\cline{2-8}
&Offline Stokes 	&10.3	&13.9	&16.9	&17.6	&19.9	&23.3\\
&Online Stokes	&2.1		&2.1		&2.3		&2.4		&2.2		&2.4\\
&Offline Navier-Stokes			&10.5	&13.6	&16.5	&17.3	&19.9	&24.0\\
&Online Navier-Stokes			&1.7		&1.8		&1.9		&2.0		&1.9		&2.0\\
\hline
\multirow{6}{*}{$128$} & \multicolumn{1} {c|}{$k$}	&237& 732&1383&2109&3039&4083  \\
	& \multicolumn{1} {c|}{$n_{deim}$}	 & 79 & 244 &461 &703 &1013 &1361 \\
\cline{2-8}
&Offline Stokes 	&8.8	&16.4	&21.0 	&17.8	&19.9	&25.4\\
&Online Stokes	&1.9	&2.2		&2.5		&2.3		&2.2		&2.4\\	
&Offline Navier-Stokes			&8.9	&16.5	&20.5	&17.9	&20.1	&25.1\\
&Online Navier-Stokes			&1.8	&1.8		&2.1		&2.1		&1.9		&2.1\\			
\hline
\end{tabular}
\end{center}
\caption{Average iteration counts of preconditioned \texttt{bicgstab} for solving 
equation (\ref{eq:deim_picard}),
over $n_s = 10$ parameter samples.}
\label{tab:it_countstab}
\end{table}


\begin{table}
\renewcommand{\baselinestretch}{1}
\begin{center}
\begin{tabular}{|c|l|l|cccccc|}
\hline
$n$& \multicolumn{2} {c}{$m$}    & 4 & 9 &16& 25 &  36 & 49 \\
\hline
\multirow{7}{*}{$32$}   & \multicolumn{2} {c|}{$k$}             &306&942        & 1485  &&& \\
                                    & \multicolumn{2} {c|}{$n_{deim}$}      &102     &314     & 495   & & & \\
\cline{2-9}
&Full   & Direct                &1.11   &1.16   &1.06   &&
&\\
&Reduced & Direct               &0.18   &1.21   &3.00   &&
&\\
&DEIM   & Direct                &\textbf{0.05}  &\textbf{0.37}  &1.02   &&
&\\
&DEIM   &Stokes &\textbf{0.05}  &0.42   &\textbf{0.95}  &&
&\\
&DEIM   & Navier-Stokes         &\textbf{0.05}  &0.44   &1.03   &&
&\\
\hline

\multirow{7}{*}{$64$} 	& \multicolumn{2} {c|}{$k$}		&273&825& 1503&2394&3339&4455 \\
					& \multicolumn{2} {c|}{$n_{deim}$}	 &91	&275& 501&798 &  1113 &1485 \\
\cline{2-9}
&Full 	& Direct		&11.0	&10.8	&10.2	&11.3	&10.1	&10.3\\
&Reduced & Direct		&0.48	&2.53	&7.33	&20.5	&39.5	&76.2\\
&DEIM  	& Direct		&\textbf{0.07}	&\textbf{0.27}	&1.08	&4.40	&9.54	&20.7\\
&DEIM 	& Stokes	&\textbf{0.07}	&0.29	&\textbf{0.86}	&2.29	&4.62	&9.04\\
&DEIM	& Navier-Stokes		&0.08	&0.30	&0.87	&\textbf{2.27}	&\textbf{4.39}	&\textbf{8.98}\\
\hline
\multirow{7}{*}{$128$} & \multicolumn{2} {c|}{$k$}	& 237&732&1383&2109 &3039&4083 \\
	& \multicolumn{2} {c|}{$n_{deim}$}	 & 79 & 244 &461 &703 &1013 &1361 \\
\cline{2-9}
&Full 	& Direct		&135&141	&147			&155&132&148\\
&Reduced & Direct		&1.62		&7.25	&23.8		&56.7&98.1&191\\
&DEIM  	& Direct		&\textbf{0.09}	&\textbf{0.39}	&1.12	&3.59&7.11&15.7\\
&DEIM 	& Stokes	&\textbf{0.09}	&0.45	&1.10		&\textbf{2.40}&\textbf{3.89}&8.74\\
&DEIM	& Navier-Stokes		&\textbf{0.09}	&0.45	&\textbf{1.08}&2.45&3.91&\textbf{8.62}\\
\hline
\end{tabular}
\end{center}
\caption{Average time for the entire nonlinear solve, over $n_s = 10$ parameter samples.} 
\label{tab:timingstab}
\end{table}	

\begin{table}
\renewcommand{\baselinestretch}{1}
\begin{center}
\begin{tabular}{|c|l|cccccc|}
\hline
$n$& \multicolumn{1} {c}{$m$}	 & 4 & 9 &16& 25 &  36 & 49 \\ 
\hline
\multirow{4}{*}{$32$}  	& \multicolumn{1} {c|}{$k$}		&306&942	& 1485 	&	&
& \\
		& \multicolumn{1} {c|}{$n_{deim}$} &102 & 314 &  495 &&& \\
\cline{2-8}
&Stokes	&0.04	&0.36	&1.02 &&
&\\
& Navier-Stokes		&1.01 & 	1.42	& 2.06	&&
&\\
\hline
\multirow{4}{*}{$64$} 	& \multicolumn{1} {c|}{$k$}		&273&825& 1503&2394&3339&4455 \\
  & \multicolumn{1} {c|}{$n_{deim}$}	 &91	&275& 501&798 &  1113 &1485 \\
\cline{2-8}
& Stokes	&0.10	&0.85	&2.38& 6.51 &13.7 &26.7	\\
& Navier-Stokes		&8.96	&10.0	&12.6&15.7 & 22.8 &35.7 \\
\hline
\multirow{4}{*}{$128$} & \multicolumn{1} {c|}{$k$}	& 237&732&1383&2109 &3039&4083 \\
& \multicolumn{1} {c|}{$n_{deim}$}	 & 79 & 244 &461 &703 &1013 &1361 \\
\cline{2-8}
& Stokes	&0.29	&	1.96 & 6.84&17.9&33.0&68.2\\
& Navier-Stokes		&116	 & 	132 & 129 &143&147&193\\
\hline
\end{tabular}
\end{center}
\caption{CPU time to construct the preconditioner}
\label{tab:offline_costs}
\end{table}
Table \ref{tab:it_countstab} illustrates that the performance of the offline preconditioners using the mean parameter compared to the versions that use the exact parameter, with the offline preconditioners requiring more iterations, as expected. In Table~\ref{tab:timingstab} we compare these offline parameters with the direct DEIM method and determine that for large enough $m$ (number of parameters) and $k$ (size of the reduced basis), the iterative methods are faster than direct methods. For example, for $m=9$, the direct methods are slightly faster whereas for $m =16$ the iterative methods are faster for all values of $n$. We also note that the fastest DEIM method is faster than the full model for all cases. Returning to the motivating example of $n=64$ and $m = 49$, the DEIM iterative method is faster than the full model, whereas the DEIM direct method performs twice as slowly as the full model. Thus, utilizing iterative methods increases the range of $k$ where reduced-order modeling is practical. 

Table \ref{tab:offline_costs} presents the (offline) cost of constructing the preconditioners. 
Since the Navier-Stokes preconditioner uses the full solution, the cost of constructing this 
preconditioner scales with the costs of the full solution. 
The cost of constructing the Stokes preconditioner is significantly smaller, because it does not require 
solving the full nonlinear problem at the mean parameter. 
It performs similarly to the exact preconditioner in the online computations.
Thus, the exact Stokes preconditioner is an efficient option for both offline and online components of this problem.

{\bf Remark.}
A variant of the offline preconditioning methods discussed here is a 
{\em blended} approach, in which a small number of preconditioners corresponding 
to a small set of parameters is constructed offline.
Then for online solution of a problem with parameter $\xi$, the preconditioner 
derived from the parameter closest to $\xi$ can be applied.
We found that this approach did not improve performance of the preconditioners 
considered here.

\section{Conclusion}
We have shown that the discrete interpolation method is effective for solving the steady-state Navier-Stokes equations. This approach produces a reduced-order model that is essentially as accurate as a naive implementation of a reduced basis method without incurring online costs of order $N$. In cases where the dimension of the reduced basis is large, performance of the DEIM is improved through the use of preconditioned iterative methods to solve the linear systems arising at each nonlinear Picard iteration. This is achieved using the mean parameter to construct preconditioners. These preconditioners are effective for preconditioning the reduced model in the entire parameter space.

\bibliography{references}{}
\end{document}